\documentclass[psamsfonts]{amsart}
\usepackage{amsmath,upref,amssymb}
\usepackage{pdfsync}
\usepackage{mathtools}
\usepackage{etoolbox}
\usepackage{scalerel}
\usepackage{float}
\usepackage{dynkin-diagrams}
\RequirePackage{calrsfs}
\usepackage[all]{xy}
\usepackage{color} 
\usepackage{eepic} 
\usepackage{multicol}
\usepackage{amsmath}
\usepackage{enumitem}
\usepackage{gensymb}
\usepackage[english]{babel}
\usepackage[T1]{fontenc}   

\usepackage{geometry}                
\usepackage{graphicx}
\usepackage{amssymb}
\usepackage{epstopdf}
\usepackage{amsmath,upref,amssymb,amsthm}
\usepackage{MnSymbol}
\usepackage{tensor}
\input xypic

\xyoption{all}

\CompileMatrices
\CompileMatrices
\numberwithin{equation}{subsection} 

\theoremstyle{plain}

\newtheorem{theorem}{Theorem}[section]

\newtheorem{corollary}[theorem]{Corollary}
   
\newtheorem{definition}[theorem]{Definition}

\newtheorem{Remark}[theorem]{Remarks}

\theoremstyle{definition}

\theoremstyle{remark}
\newtheorem*{note}{Note}

\makeatletter
\renewcommand{\author}[2][]{%
  \def\@tempa{#1}
  \ifx\@empty\authors
    \ifx\@tempa\@empty
      \gdef\shortauthors{#2}%
    \else
      \gdef\shortauthors{#1}%
    \fi
    \gdef\authors{\author{#2}}%
  \else
    \ifx\@tempa\@empty
      \g@addto@macro\shortauthors{\and#2}%
    \else
      \g@addto@macro\shortauthors{\and#1}%
    \fi
    \g@addto@macro\authors{\and\author{#2}}%
  \fi
}
\renewcommand{\address}[2][]{\g@addto@macro\authors{\address{#1}{#2}}}
\def\@setauthors{%
  \begin{center}%
    \footnotesize
    \vspace{20pt}
    \let\and\@empty
    \def\author##1{\advance\@tempcnta\@ne}%
    \def\address##1##2{\advance\@tempcntb\@ne}%
    \@tempcnta=\z@  \@tempcntb=\z@
    \authors
    \ifnum\@tempcnta>\@ne \ifnum\@tempcntb=\@ne
        \oneaddress
      \else
        \sepaddresses
      \fi
    \else
      \oneaddress
    \fi
  \end{center}%
}
\def\oneaddress{%
  \begingroup
  \let\author\@iden \let\address\@gobbletwo
  \renewcommand{\andify}{%
    \nxandlist{\unskip, }{\unskip{} and~}{\unskip, and~}}%
  \uppercasenonmath\authors
  \andify\authors
  \authors
  \endgroup
  \begingroup \let\and\relax \let\author\@gobble
  \def\address##1##2{\unskip\\[10pt] \itshape##2}%
  \authors
  \endgroup
}
\def\sepaddresses{%
  \begingroup
    \baselineskip10\p@\relax
    \def\address##1##2{ ({\itshape##2}\/)}
    \def\author##1{\def\temp{##1}\leavevmode\uppercasenonmath\temp\temp}%
    \nxandlist
      {,\\[\baselineskip]}
      {\\[\baselineskip] \textsc{\lowercase{and}}\\[\baselineskip]}
      {,\\[\baselineskip]\textsc{\lowercase{and}}\\[\baselineskip]}
      \authors 
    \authors
  \endgroup
}
\def\maketitle{\par
  \@topnum\z@
  \@setcopyright
  \thispagestyle{firstpage}%
  \uppercasenonmath\shorttitle
  \ifx\@empty\shortauthors \let\shortauthors\shorttitle
  \else
    \newcommand{\@xuppercasenonmath}[1]{\toks@\@emptytoks
      \@xp\@skipmath\@xp\@empty##1$$%
      \edef##1{\@nx\protect\@nx\@upprep\the\toks@}}%
    \@xuppercasenonmath\shortauthors
    \def\@@and{AND}
    \renewcommand{\andify}{%
      \nxandlist{\unskip, }{\unskip{ }\@@and{ }}{\unskip, \@@and{ }}}%
    \andify\shortauthors
  \fi
  \@maketitle@hook
  \begingroup
  \@maketitle
  \endgroup
  \c@footnote\z@
  \@cleartopmattertags
}
\def\@maketitle{%
  \normalfont\normalsize
  \let\@makefntext\noindent
  \@adminfootnotes
  \ifx\@empty\addresses\else \@footnotetext{\@setotheraddresses}\fi
  \global\topskip68\p@\relax
  \@settitle
  \ifx\@empty\authors \else \@setauthors \fi
  \ifx\@empty\@dedicatory
  \else
    \baselineskip26\p@
    \vtop{\centering{\footnotesize\itshape\@dedicatory\@@par}%
      \global\dimen@i\prevdepth}\prevdepth\dimen@i
  \fi
  \toks@\@xp{\shortauthors}\@temptokena\@xp{\shorttitle}%
  \edef\@tempa{\@nx\markboth{\the\toks@}{\the\@temptokena}}\@tempa
  \@setabstract
  \normalsize
  \if@titlepage
    \newpage
  \else
    \dimen@34\p@ \advance\dimen@-\baselineskip
    \vskip\dimen@\relax
  \fi
} 
\renewcommand{\thanks}[1]{%
  \ifx\@empty\thankses
    \gdef\thankses{\thanks{#1}}%
  \else
    \g@addto@macro\thankses{\endgraf\thanks{#1}}%
  \fi}
\def\@setthanks{\def\thanks##1{\noindent##1\@addpunct.}\thankses}
\renewcommand{\curraddr}[2][]{%
  \ifx\@empty\addresses
    \gdef\addresses{\curraddr{#1}{#2}}%
  \else
    \g@addto@macro\addresses{\endgraf\curraddr{#1}{#2}}%
  \fi}
\renewcommand{\email}[2][]{%
  \ifx\@empty\addresses
    \gdef\addresses{\email{#1}{#2}}%
  \else
    \g@addto@macro\addresses{\endgraf\email{#1}{#2}}%
  \fi}
\renewcommand{\urladdr}[2][]{%
  \ifx\@empty\addresses
    \gdef\addresses{\urladdr{#1}{#2}}%
  \else
    \g@addto@macro\addresses{\endgraf\urladdr{#1}{#2}}%
  \fi}
\def\@setotheraddresses{%
  \def\curraddr##1##2{\noindent
    \emph{Current address\@ifnotempty{##1}{ of ##1}}:\space
      ##2\@addpunct.}%
  \def\email##1##2{\noindent
    \emph{E-mail address\@ifnotempty{##1}{ of ##1}}:\space
      \texttt{##2}}%
  \def\urladdr##1##2{\noindent
    \emph{WWW address\@ifnotempty{##1}{ of ##1}}:\space
      \texttt{##2}}%
  \addresses
}
\let\enddoc@text\relax
\makeatother

\begin{document}

\title{Subregular J-Rings of the Finite Irreducible Coxeter Systems }

\author {Annette Pilkington}
 \address{Department of Mathematics,  University of Notre
Dame,  Notre Dame, Indiana, 46556}
\email{Pilkington.4@nd.edu}

\maketitle    

\begin{abstract}
In  previous papers, the author showed that in many cases of interest there exists  an isomorphism between certain  path algebras  related to  the structure of the subregular J-rings of Coxeter systems 
 and  matrix rings
over a free product of rings.  For the finite irreducible Coxeter systems, the application of the  isomorphism theorems was straightforward in all but a few cases. 
 In  this paper we show how the isomorphism theorems can be applied to show that 
 the path algebras in the remaining cases are isomorphic to a direct sum of matrix rings over 
 the rational numbers or over a direct product of field extensions of the rational numbers.
\end{abstract}

\section{Introduction}
In Pilkington \cite{Me1} and \cite{Me2}, we showed that  in a wide range of cases certain path algebras, $\widehat{R}_{W, S}$ and 
$\widehat{\widetilde{R}}_{W, S}$(described below), associated to 
a Coxeter system, $(W, S)$, are isomorphic to matrix rings over the free product of rings with generators which correspond to the edges of the Coxeter graph and the corresponding entries in the Coxeter matrix. The algebra $\widehat{\widetilde{R}}_{W, S}$ is strongly related to the  subregular J-ring of the Coxeter system $(W, S)$. The application of the theorems introduced in Pilkington \cite{Me2} to the structure of these path algebras for the Coxeter systems associated to the finite irreducible 
Coxeter groups was straightforward for all but three cases, namely that of the algebra 
$\widehat{\widetilde{R}}_{W, S}$ for 
  Coxeter systems, $(W, S)$ of type $I_2(n), n > 5, n$ even,  $B_n$, $n \geq 2$, and $F_4$. In Pilkington \cite{Me2}, Example 5.9,  we showed that in the case of $B_4$ there was a unital ring homomorphism between the algebra $\widehat{\widetilde{R}}_{W, S}$ and  
  an algebra to which the theorems did apply. In this paper we show that a similar unital ring homomorphism exists more generally and can be applied to $\widehat{\widetilde{R}}_{W, S}$
  for the exceptional cases listed above to show that $\widehat{\widetilde{R}}_{W, S}$ is isomorphic to a direct sum of matrix algebras over $\Bbb{Q}$ or over a direct product of field extensions of $\Bbb{Q}$ in these cases.

Let $(W, S= \{s_i\}_{1 \leq i \leq N})$ be a Coxeter system of finite rank with Coxeter matrix $(m_{ij})_{1 \leq i, j \leq N}$. The  graph $\Gamma_{W, S} =(S, Y)$ associated to the Coxeter system $(W, S)$  has vertices $S$ and  edges $Y = \{y_{ij} = [s_is_j] \  |\  1 \leq i, j \leq N, \ 3 \leq m_{ij}  \leq \infty\}$. We let  $P_{W,S}$ denote the path algebra over $\mathbb{Z}$  associated to the  graph $\Gamma_{W,S}$, and we let $\widehat{P}_{(W, S)} = P_{W,S}\otimes \Bbb{Q}$.

For $n \geq 3$, we denote  the minimal polynomial of $\displaystyle 4\cos^2\frac{\pi}{n}$ over $\Bbb{Q}$ 
by $C_n(t) \in \Bbb{Z}[t]$ and we let $$c_n(t) = \prod_{{\tiny \left.\begin{array}{c}N\in\mathbb{N}_{\geq 3} \\N\mid n\end{array}\right.}}C_N(t).$$
 For example, we have 
$$\{C_3(t) = t - 1, \ \ C_4(t) = t - 2, \ \ C_5(t) = t^2 - 3t + 1, \ \ C_6(t) = t - 3, \dots \}.$$
and
$$\{c_3(t) = t - 1, \  \ c_4(t) = t - 2, \ \ c_5(t) = t^2 - 3t + 1, \ \ c_6(t) = (t - 1)(t - 3) = t^2 - 4t + 3, \dots \}.$$
We let  $I_{W,S}$ be the two-sided ideal of $P_{W,S}$ generated by the set 
$$\left\{C_{m_{ij}}([s_js_is_j]) = C_{m_{ij}}(y_{ji}y_{ij}) \ |\ (s_i,s_j)\in S\times S, \ 2 < m_{ij} < \infty\right\},$$ 
where $f(y_{ji}y_{ij})$ denotes the evaluation of the polynomial $f(t)$ at $y_{ji}y_{ij}$  in the ring 
$[s_j]P_{(W, S)}[s_j]$ (this ring has identity element $[s_j]$). 
We let $\widetilde{I}_{W,S}$ denote the two-sided ideal of $P_{W,S}$ generated by the elements of the union of the sets 
$\left\{c_{m_{ij}}(y_{ji}y_{ij}) \ |\ (s_i,s_j)\in S\times S, \ 2 < m_{ij} < \infty,
m_{ij} \ \mbox{odd}\right\}$
and $ \left\{y_{ij}c_{m_{ij}}(y_{ji}y_{ij}) \ |\ (s_i,s_j)\in S\times S, \ 2 < m_{ij} < \infty,
m_{ij} \ \mbox{even}\right\}$. 
Define the quotient algebras $R_{W,S}:=P_{W,S}/{I}_{W,S}$ and 
$\widetilde{R}_{W,S}:=P_{W,S}/\widetilde{I}_{W,S}$ of $P_{W,S}$.
We let $\widehat{R}_{W,S} = R_{W,S}\otimes \Bbb{Q}$ and 
$\widehat{\widetilde{R}}_{W,S} = \widetilde{R}_{W,S}\otimes \Bbb{Q}$.

 The latter algebra is of considerable interest in the representation theory of Coxeter groups. The algebra $\widetilde{R}_{W, S}$ is essentially the subregular $J$-ring of a Coxeter system, that is, the two-sided ideal $J_1$ of Lusztig's asymptotic Hecke  algebra spanned by its standard basis elements parameterized by elements of {\bf a}-value $1$ in the Coxeter group (see \cite{Unequal}, \cite{Noncom} and  \cite{Subregular} for precise statements).  In the case where all finite entries not equal to 1 of the Coxeter matrix are prime, both algebras, $R_{W, S}$ and 
$\widetilde{R}_{W, S}$,  coincide.

In  Pilkington \cite{Me1} and \cite{Me2}, we showed that in situations where the  Coxeter graph has  suitable maximal trees,  the algebras $\widehat{R}_{W, S}$  and $\widehat{\widetilde{R}}_{W', S'}$
are isomorphic to the $N \times N$ matrices  over a free product of rings related to the edges of  the Coxeter graph, where $N = |S|$. These rings are either field extensions of $\Bbb{Q}$, rings of Laurent polynomials over $\Bbb{Q}$ or rings of polynomials in two non-commuting variables over $\Bbb{Q}$. They may also  include direct products of field extensions of $\Bbb{Q}$ or   quotient rings of 
 rings of polynomials in two non-commuting variables over $\Bbb{Q}$ in the case of $\widehat{\widetilde{R}}_{W, S}$. In Pilkington \cite{Me2}, we applied these results to the finite irreducible Coxeter groups, with the following result:\\
 
 \noindent
{\bf Theorem} ( {\bf Pilkington \cite{Me2}, Theorem 5.8})
 Let $(W, S)$ be a Coxeter system with rank $N$, where $W$ is a finite irreducible Coxeter group.  We have the following unital ring isomomorphisms:
  $$\widehat{R}_{W,S} \cong \begin{cases} M_N(\Bbb{Q}) & \mbox{if $(W, S)$ is of type $A_n (n \geq 1), B_n (n \geq 2), D_n (n \geq 4), E_6, E_7, E_8, F_4$ or $G_2$,}\\
M_N(\Bbb{Q}_5) & \mbox{if $(W, S)$ is of type $H_3,$ or $H_4$,}\\
M_2(\Bbb{Q}_n) & \mbox{if $(W, S)$ is of type $I_2(n), n \geq 5, n \not= 6$,}\\
\end{cases}$$
(where $\Bbb{Q}_n$ if the field extension of $\Bbb{Q}$ obtained by adjoining $4\cos^2\left(\frac{\pi}{n}\right)$) \\
and
 $$\widehat{\widetilde{R}}_{W,S} \cong \begin{cases} M_N(\Bbb{Q}) & \mbox{if $(W, S)$ is of type $A_n (n \geq 1), D_n (n \geq 4), E_6, E_7$ or $E_8$,}\\
M_N(\Bbb{Q}_5) & \mbox{if $(W, S)$ is of type $H_3,$ or $H_4$,}\\
M_2({L}_n) & \mbox{if $(W, S)$ is of type $I_2(n), n \geq 5$, where $n$ is odd,}\\
\end{cases}$$
(where, for $n \geq 3$, $L_n$ denotes the  following  direct product of field extensions of $\Bbb{Q}$; 
$$L_n = \prod_{{\tiny \left.\begin{array}{c}N\in\mathbb{N}_{\geq 3} \\N\mid n\end{array}\right.}}\Bbb{Q}_N).$$
In Pilkington \cite{Me2}, Example 5.9,  we showed  that if  $(W, S)$ is a Coxeter system of type $B_4$, we have the following unital ring  isomorphism:
$$\widehat{\widetilde{R}}_{W,S} \cong M_4(\Bbb{Q}) \bigoplus M_{3}(\Bbb{Q}) \bigoplus \Bbb{Q}.$$
In this paper, we show how  techniques similar to those used in  Pilkington \cite{Me2}, Example 5.9 apply more generally. We show that :\\

\noindent
{\bf Theorem} \ref{general}
Let  $(W, S)$ be a Coxeter system  where $S = \left\{s_i, t_j \ \vert \ 1 \leq i \leq n + 1, 1 \leq j \leq m + 1\right\}$, with graph $\Gamma_{W, S}$ shown below, where $a \geq 4$ is an even integer, and  the labels on the edges correspond to the associated entries of the Coxeter matrix. 

\begin{equation*}
\xymatrix{
s_1\ar@{-}[r]^{3}&s_2\ar@{-}[r]^{3}&{s_3}\ar@{.}[r]&{s_{n}}\ar@{-}[r]^{3}&{s_{n + 1}}\ar@{-}[r]^{a}&{t_{m + 1}}\ar@{-}[r]^{3}&t_{m}\ar@{.}[r]&{t_{3}}\ar@{-}[r]^{3}&{t_2}\ar@{-}[r]^{3}&{t_1}\\
}
\end{equation*}
 \noindent
Then we have a unital ring isomorphism; 
$$\widehat{\widetilde{R}}_{W,S}  \cong 
  M_{m + n + 2}({L}_a) \bigoplus M_{n}(\Bbb{Q}) \bigoplus M_{m}(\Bbb{Q}).$$
\\

Applying the above theorem   to the algebra  $\widehat{\widetilde{R}}_{W,S} $ 
for the Coxeter systems of type  $G_2 = I_2(6)$, $I_2(n), n > 6, n$ even,  $B_n, n \geq 2$ and $F_4$, the cases missing from Pilkington \cite{Me2}, Theorem 5.8, we  extend the theorem as follows:\\

  \noindent
{\bf Theorem} \ref{complete}
 Let $(W, S)$ be a Coxeter system with rank $N$, where $W$ is a finite irreducible Coxeter group (the Coxeter graphs of such systems are shown in Figure  \ref{fig:fig1}). We have the following unital ring isomomorphisms:
 $$\widehat{R}_{W,S} \cong \begin{cases} M_N(\Bbb{Q}) & \mbox{if $(W, S)$ is of type $A_n (n \geq 1), B_n (n \geq 2), D_n (n \geq 4), E_6, E_7, E_8, F_4$ or $G_2$,}\\
M_N(\Bbb{Q}_5) & \mbox{if $(W, S)$ is of type $H_3,$ or $H_4$,}\\
M_2(\Bbb{Q}_n) & \mbox{if $(W, S)$ is of type $I_2(n), n \geq 5, n \not= 6$,}\\
\end{cases}$$
\\
and
 $$\widehat{\widetilde{R}}_{W,S} \cong \begin{cases} M_N(\Bbb{Q}) & \mbox{if $(W, S)$ is of type $A_n (n \geq 1), D_n (n \geq 4), E_6, E_7$ or $E_8$,}\\
M_N(\Bbb{Q}_5) & \mbox{if $(W, S)$ is of type $H_3,$ or $H_4$,}\\
M_2({L}_n) & \mbox{if $(W, S)$ is of type $I_2(n), n \geq 5$, where $n$ is odd,}\\
M_2({L}_n) \oplus \Bbb{Q} \oplus \Bbb{Q} & \mbox{if $(W, S)$ is of type $I_2(n), n \geq 6$, where $n$ is even,}\\
M_n(\Bbb{Q}) \oplus M_{n-1}(\Bbb{Q}) \oplus \Bbb{Q} & \mbox{if $(W, S)$ is of type $B_n$, $n \geq 2$,}\\
M_4(\Bbb{Q}) \oplus M_{2}(\Bbb{Q}) \oplus M_{2}(\Bbb{Q}) & \mbox{if $(W, S)$ is of type $F_4$,}\\
\end{cases}$$
where $\Bbb{Q}_n$ and  $L_n$  are as described above. \\

In Section 2, we review the statement of some necessary definitions and theorems from Pilkington 
\cite{Me2} and in Section 3, we look at the applications to   Coxeter systems.

\section{Graphs, Maximal Trees and Path Algebras}
As in Pilkington \cite{Me1}, \cite{Me2}, we use  the notation and terminology of Serre \cite{Serre} for graphs and trees. This somewhat  unconventional definition can be thought of as an undirected graph 
with  labels for the edges in both directions.  
\begin{definition}\label{defgraph}
A graph  $\Gamma = (S, Y)$ consists of a set of vertices $S = \mbox{vert} \ \Gamma $, a set of edges  $Y = \mbox{edge} \ \Gamma$ and two maps 
$$Y \to S \times S, \ \  y \to (o(y), t(y))$$
and 
$$Y \to Y, \ \ y \to \bar{y},$$
which satisfy the following condition: for each $y \in Y$, we have $\bar{\bar{y}} = y, \bar{y} \not= y$ and 
$o(y) = t(\bar{y})$. 
\end{definition}

Each $s \in S$ is called a {\it vertex} of $\Gamma$ and each $y \in Y$ is called an {\it edge} of $\Gamma$. 
If $y \in Y$, $o(y)$ is called the {\it origin} of $y$ and $t(y)$ is called the {\it terminus} of $y$, together $o(y)$ and $t(y)$ are called the {\it extremities }
of $y$. 
We will represent such a graph by a diagram where a point on the diagram represents a vertex and a single line joining two points represents a set of edges of the form $\{y, \bar{y}\}$. 

A graph is said to have {\it multiple edges} if two edges have the same origin and terminus. 
A {\it loop} is an edge with the same origin and terminus. 
 A graph with no loops and no multiple edges is called a {\it combinatorial graph}. 
\begin{definition}\label{orientdef} Let $\Gamma = (S, Y)$ be a graph. An {\it orientation} of  $\Gamma$ is a subset $Y_+$ of $Y$  such that $Y$ is the disjoint union 
$Y = Y_+ \bigcupdot \overline{Y_+}$. 
\end{definition}
\subsection{Paths}
A {\it path of length $n \geq 1$} in a graph $\Gamma = (S, Y)$, where $S = \{s_i\}_{i \in I}$ is indexed by the set $I$,  is a sequence of edges $p = (y_1, \dots , y_n)$ such that $t(y_i) = o(y_{i + 1})$ for   $1 \leq i \leq n - 1$. 
We will  denote such a path in what follows as $p = y_1 \dots y_n$. {\it Paths of length zero}  are just single vertices and will be denoted by $p = [s]$ for $s \in S$.
 
 Let 
$\frak{P}$ denote the set of all paths in $\Gamma$, 
we can extend the maps $o$ and $t$ to $\frak{P}$ by letting $o( y_1 \dots  y_n ) = o(y_1) $ and $t( y_1 \dots  y_n ) = t(y_n)$ and $o([s]) = s = t([s])$. 
 If $s_{i_k} = t(y_k)$, $s_{i_{k - 1}} = o(y_k)$, where $i_k \in I$, for $0 \leq k \leq n$, we say that $p =  y_1 \dots  y_n$ is a path from $s_{i_0}$ to $s_{i_n}$ and that  $s_{i_0}$ and  $s_{i_n}$ are the {\it extremities} of the path. 
    If the graph $\Gamma$ is combinatorial, a path $p = y_1\dots y_n$ is determined by the extremities of its edges and can be characterized and denoted by its vertices as 
 $$p = y_1 \dots y_n = [s_{i_0}s_{i_1} \dots s_{i_n}], \ \mbox{where} \ s_{i_k} = t(y_k), s_{i_{k - 1}} = o(y_k),\ \mbox{for} \ i_k \in I, 0\leq k \leq n.$$  Coxeter graphs are combinatorial and 
 we will use both sets of notation to describe paths in the associated algebras. 
The  graph $\Gamma$ is {\it connected} if any two vertices are the extremities of at least one path.
 If $p = y_1 \dots y_n$ and $y_{i + 1} = \overline{y_i}$
for some $1 \leq i \leq n - 1$, then the pair $(y_i, y_{i + 1})$ is called a {\it backtracking}.  A path $p$ is a {\it circuit} if it is a path 
 without backtracking such that $o(p) = t(p)$ and $\ell(p) \geq 1$.

 \subsection{Trees}\label{trees} If $T \not= \emptyset$ is a connected graph without circuits, $T$ is called a {\it tree}.

\begin{definition} \label{subgraph} A subgraph of a graph $\Gamma = (S, Y)$ is a graph $\Gamma' = (S', Y')$ 
such that $S' \subseteq S$ and $Y' \subseteq Y$, where $\{ o(y), t(y) \ |\  y \in Y'\} \subseteq S'$,
 and $\overline{Y'} = Y'$.  The maps $Y' \to S' \times S', y \to (o(y), t(y))$ and $Y' \to Y', y \to \bar{y}$
are inherited from $\Gamma$. 
\end{definition}

Given a non-empty connected graph $\Gamma$, the set of subgraphs  which are trees are partially ordered by inclusion.  The union of a chain of trees is a tree and  hence the set has a maximal element, (using   Zorn's lemma if infinite). Such a  maximal element is called a {\em maximal tree } of $\Gamma$. 

\subsection{Path Algebras}\label{palgebras}
Let $\Gamma = (S, Y)$ be a graph with vertices  $S = \{s_i\}_{i \in I}$, edges $Y$ and paths $\mathfrak {P}$. 
 Given a commutative ring $A$ with identity, the {\it path algebra of $\Gamma$ over $A$}, denoted $A\Gamma$, is an associative $A$-algebra, with an $A$-basis  of paths in $\mathfrak {P}$. Multiplication of paths is given by concatenation;
 $$(y_1y_2\dots y_n)(y'_1y'_2\dots y'_m) = 
 \begin{cases} 
  0 & \mbox{if} \ t(y_n)  \not= o(y'_1)\\ 
  y_1y_2\dots y_ny'_1y'_2\dots y'_m & \mbox{if} \ t(y_n)  = o(y'_1)\\ 
  \end{cases}$$

Under this multiplication rule, the paths of length $0$, $\{[s_i]\}_{s_i \in S}$,  form a set of  orthogonal idempotents. If $S$ is finite,  then $A\Gamma$ has an identity given by 
$$1 = \sum_{s \in S}[s].$$

In what follows, if $f \in A[t]$ and $y \in Y$, we let $f(y\bar{y})$ denote the element of the subring $[o(y)]A\Gamma[o(y)] \subset A\Gamma$  (where the imbedding is non-unital)
that we get by  replacing $t$ by $y\bar{y}$ and $1_A$ by $[o(y)]$ in the polynomial $f$ (this is $f$ evaluated at $y\bar{y}$ in the subring $[o(y)]A\Gamma[o(y)]$). 

\begin{definition}\label{propF} Let  $Y_1$ be a subset of $Y$ such that $\overline{Y_1} = Y_1$. We say that a set of polynomials, $\mathfrak{F}_{Y_1} = \{f_y(t)\}_{y \in Y_1 }$,  in the polynomial ring $A[t]$ has {\it property F} if the following  conditions are satisfied for each $y \in Y_1$:
\begin{itemize}
\item $f_y(t) = th_y(t) - \kappa_y $,  where $h_y(t)$
 is a non-zero polynomial in $A[t]$ and $\kappa_y$ is a unit in $A$. 
\item $f_y = f_{\bar{y}}$.
\end{itemize}
\end{definition}

Now consider a subset $Y_1$  of $Y$, such that  $\overline{Y_1} = Y_1$,  with associated set of polynomials 
$\mathfrak{F}_{Y_1} =  \{f_y(t)\}_{y \in Y_1 }$, where $\mathfrak{F}_{ Y_1}$ has property $F$. 
We let $I_{{\mathfrak{F}_{Y_1}}}$ denote the two-sided ideal 
of $A\Gamma$ generated by the set of elements $\{f_y(y\bar{y}) \ |\  y \in Y_1\}$. 
We let $R_{{\mathfrak{F}_{Y_1}}}$ denote the quotient algebra
$$R_{{\mathfrak{F}_{Y_1}}} = A\Gamma/I_{{\mathfrak{F}_{Y_1}}}.$$

\begin{Remark} \label{relations}\ 

\begin{itemize}
\item Note that each pair $\{y, \bar{y}\} \subset Y_1$ gives two relations in $R_{\mathfrak{F}_{Y_1}}$,  namely $f_y(y\bar{y}) = 0$ and $f_{\bar{y}}(\bar{y}y) = 
f_y(\bar{y}y)= 0$.
\item If the graph has a single vertex, $[s]$, and no edges, that is $\Gamma = (S = \{s\}, Y = \emptyset)$, 
then $Y_1  = \emptyset$ and the set of generators for the ideal $ I_{{\mathfrak{F}_{Y_1}}}$ is empty. In this case there is a unital ring isomorphism; 
 $R_{{\mathfrak{F}_{Y_1}}}  = A[s] \cong A$.  
\item If the graph, $\Gamma = (S, Y)$, has connected components $\Gamma^j = (S^j, Y^j),  j \in J$, 
then $\displaystyle R_{{\mathfrak{F}_{Y_1}}}  \cong \bigoplus_{j \in J}R^j$, where 
$R^j = A\Gamma^j/I_{{\mathfrak{F}_{Y_1\cap Y^j}}}$. This follows from the observations that $\displaystyle A\Gamma = \bigoplus_{j \in J}A\Gamma^j$ and 
$I_{{\mathfrak{F}_{Y_1}}}$ is the direct sum of the ideals
$I_{{\mathfrak{F}_{Y_1\cap Y^j}}}$ for  $j \in J$. 
\end{itemize}
\end{Remark}

\subsection{Free Products of Rings}
\begin{definition}\label{definition} (Cohn \cite{Cohn1}, Section 3)
Let  $K$ be a field. A  $K$-algebra $R$ is said to be the free product of the $K$-algebras $\{R_i\}_{i \in I}$ (over $K$)  if 
$\{R_i\}_{i \in I}$ forms a family of subrings of $R$
 such that 
\begin{enumerate}[label=(\roman*)]
\item $R_i \cap R_j = K$ for $i \not= j$, 
\item if $X_i$ is a set of generators of $R_i$, $i \in I$, then $\bigcup_IX_i$ is a set of generators of $R$,
\item If $C_i$ is a set of defining relations of $R_i, i \in I$ (in terms of the set of generators $X_i$) then $\bigcup_IC_i$ is a set of defining relations of $R$ (in terms of the set of generators $\bigcup_IX_i$). 
\end{enumerate}
We denote the free product of $R_1$ and $R_2$ over $K$ by $R_1\ast_KR_2$ and the free product of the $K$-algebras $\{R_i\}_{i \in I}$ over $K$ by $(\ast_KR_i)_{i \in I}$.
\end{definition}

\begin{note}The free product of algebras over commutative rings does not always exist. However, it is guaranteed to exist when the base ring is a field, 
see Cohn \cite{Cohn1}, Corollary to Theorem 4.7. 
 When it exists, it is isomorphic to the universal product of the rings.
\end{note}

When the graph $\Gamma$ has a suitable maximal tree, we have the following isomorphism of $K$-algebras:

\begin{theorem}\label{theisocor}(Pilkington \cite{Me2}, Theorem 4.15, Lemma 4.17  and Corollary 4.18)
 Let $K$ be a field and let $\Gamma = (S, Y)$ be a connected  graph with a finite number of vertices and edges, and orientation $Y_+$, with 
  $Y = Y_+ \bigcupdot \overline{Y_+}$. Let $Y_1 \subseteq Y$ with  $\overline{Y_1} = Y_1$.
    Let $\mathfrak{F}_{Y_1} = \{f_{y}(t)\}_{y \in Y_1}$ be a set of polynomials in $K[t]$ with property F.
  Let  $T = (S, Y^T)$ be a maximal tree for $\Gamma$,  with edges $Y^T \subseteq Y_1$, and  orientation $Y^T_+ = Y_+ \cap Y^T$.
 Then $R_{\mathfrak{F}_{Y_1}}$ is isomorphic to the matrix ring $M_N(Q)$, where $N$ is the cardinality of $S$, and  $Q$ is the free product of the following rings:
 \begin{equation*}
 \begin{cases}
 K[t]/\langle f_{y}(t)\rangle & y \in Y_1\cap Y_+\\
  K[z_y, {z_y}^{-1}]& y \in (Y_1\cap Y_+)\backslash Y^T_+\\
  K\langle  x_y, \overline{x_y}\rangle& y \in Y_+ \backslash Y_1.\\
  \end{cases}
 \end{equation*}   
\end{theorem}

\section{Coxeter Systems and Associated Path Algebras}\label{coxsys}
In this section, we look at the applications of Theorem 2.7 to finite rank Coxeter Systems. In particular, we look at the applications to  the finite irreducible Coxeter systems omitted from Pilkington \cite{Me2}, Theorem 5.8,  namely $G_2, I_2(n), n > 6, n \ \mbox{even}, B_n, n \geq 2$ and $F_4$.

\subsection{Coxeter Systems}
A finite rank 
  {\em Coxeter matrix} is a matrix $M = (m_{ij})_{1 \leq i, j \leq N}$,  with 
$m_{ii} = 1$  and $m_{ij} = m_{ji} \in \mathbb{N}_{\geq 2}\cup \{\infty\}$ for all  $1 \leq  i, j  \leq N$ such that  $i \not= j$. 
For  $S = \{s_i\}_{1 \leq i \leq N}$, the corresponding {\em Coxeter group} is the group with 
presentation
$$W = \langle s_i (1\leq i \leq N) \ |\  (s_ks_l)^{m_{kl}} = 1 \ \mbox{for} \ s_k, s_l \in S \ \mbox{with} \ m_{kl} \not= \infty \rangle.$$
The pair $(W, S)$ is called a {\em Coxeter system} with Coxeter matrix $M$. 

We can associate to a  Coxeter system, $(W, S)$, a graph $\Gamma_{W, S} = (S, Y)$ with  vertices $s_i \in S$  and 
 edges  $Y = \left\{y_{ij}  \ |\ \  i \not= j \ \mbox{and} \  m_{ij} \geq 3\right\}$.
 We have $o(y_{ij}) = [s_i]$, $t(y_{ij}) = [s_j]$, and $\overline{y_{ij}} = y_{ji}$,
 for all  $1 \leq  i, j  \leq N$ such that  $i \not= j$.  The {\it Coxeter graph} of a  Coxeter system, $(W, S)$, has a point for each vertex and a single line representing each edge in the set $Y$ above, with a label $m_{ij}$ for each such edge with 
 $m_{ij} \geq 4$, that is, $m_{ij} = 3$ for the unlabelled edges.

  In Figure \ref{fig:fig1} below, we show the  Coxeter graphs of the Coxeter Systems associated  to the finite irreducible Coxeter groups, with those associated to  the Weyl groups on the left.
The graphs associated to the Coxeter Systems of type   $A_n, B_n$ and $D_n$ each have $n$ vertices.
\begin{figure}[H]
$$
\left.\begin{array}{ll}
A_n \ (n \geq 1) & \dynkin[Coxeter]A{} \\
B_n \ (n \geq 2) & \dynkin[Coxeter]B{} \\
D_n \ (n \geq 4)& \dynkin[Coxeter]D{}\\
E_6 & \dynkin[Coxeter]E6\\
E_7 & \dynkin[Coxeter]E7\\
E_8 & \dynkin[Coxeter]E8\\
F_4 & \dynkin[Coxeter]F4\\
G_2 & \dynkin[Coxeter, gonality=6]G2\\
\end{array}\right.
\qquad
\left.\begin{array}{ll}
H_3 & \dynkin[Coxeter]H3\\
H_4 & \dynkin[Coxeter]H4\\
I_2(n) \  (n \geq 5, n \not= 6) & \dynkin[Coxeter, gonality=n]I{}\\
&\\
&\\
&\\
&\\
&\\
&\\
&\\
\end{array}\right.
$$
\caption{Coxeter Graphs of the Finite Irreducible Coxeter Groups}
\label{fig:fig1}
\end{figure}

 We let $\mathfrak {P}_{W, S}$ denote the set of paths in $\Gamma_{W, S}$.
 Since the graph is combinatorial, one can describe edges and paths in the graph uniquely 
 by the sequence of vertices that they pass through. We will also frequently use $y_{ij}$ to denote the edge $[s_is_j] \in Y$. 
 We let ${P}_{W, S}$ denote the path algebra $\mathbb{Z}\Gamma_{W, S}$ and we let 
 $\widehat{P}_{W, S}$ denote the path algebra $\mathbb{Q}\Gamma_{W, S} \cong {P}_{W, S}\otimes \Bbb{Q} $.

 \subsection{The rings $ R_{W,S}$ and $ \widetilde{R}_{W,S}$.}\label{noncomrep}

\begin{definition}
For $n \geq 3$, we denote  the minimal polynomial of $\displaystyle 4\cos^2\frac{\pi}{n}$ over $\Bbb{Q}$ 
by $C_n(t) \in \Bbb{Z}[t]$ and we let $$c_n(t) = \prod_{{\tiny \left.\begin{array}{c}N\in\mathbb{N}_{\geq 3} \\N\mid n\end{array}\right.}}C_N(t).$$
\end{definition} 
\begin{Remark}\label{seen}
For all  $n \geq 3$,
 the polynomial $C_n(t)$ is of the form $tg_n(t) - \gamma_n$ for some 
$g_n(t) \in \Bbb{Z}[t]$ and some $\gamma_n \in \Bbb[Z]$ and 
 the polynomial $c_n(t)$ is of the form $th_n(t) - \kappa_n$ for some 
$h_n(t) \in \Bbb{Z}[t]$ and some $\kappa_n \in \Bbb[Z]$. (See Pilkington \cite{Me2}, Remark 3.11).
\end{Remark}

We let  $I_{W,S}$ be the two-sided ideal of $P_{W,S}$ generated by the set 
$$\left\{C_{m_{ij}}([s_js_is_j]) = C_{m_{ij}}(y_{ji}y_{ij}) \ |\ (s_i,s_j)\in S\times S, \ 2 < m_{ij} < \infty\right\},$$ 
where $f(y_{ji}y_{ij})$ denotes the evaluation of the polynomial $f(t)$ at $y_{ji}y_{ij}$  in the ring 
$[s_j]P_{W, S}[s_j]$ (this ring has identity element $[s_j]$). 
We let $\widetilde{I}_{W,S}$ denote the two-sided ideal of $P_{W,S}$ generated by the elements of the union of the sets 
$\left\{c_{m_{ij}}(y_{ji}y_{ij}) \ |\ (s_i,s_j)\in S\times S, \ 2 < m_{ij} < \infty,
m_{ij} \ \mbox{odd}\right\}$
and $ \left\{y_{ij}c_{m_{ij}}(y_{ji}y_{ij}) \ |\ (s_i,s_j)\in S\times S, \ 2 < m_{ij} < \infty,
m_{ij} \ \mbox{even}\right\}$. 
We
define the quotient algebras $R_{W,S}:=P_{W,S}/{I}_{W,S}$ and 
$\widetilde{R}_{W,S}:=P_{W,S}/\widetilde{I}_{W,S}$ of $P_{W,S}$.
We let $$\widehat{R}_{W,S} = R_{W,S}\otimes \Bbb{Q} \  \mbox{and} \ 
\widehat{\widetilde{R}}_{W,S} = \widetilde{R}_{W,S}\otimes \Bbb{Q}.$$

\begin{definition}
For $n \geq 3$, we  let $\Bbb{Q}_n = \Bbb{Q}[t]/\langle C_n(t)\rangle$ denote the field extension $\displaystyle \Bbb{Q}\left[4\cos^2\left(\frac{\pi}{n}\right)\right]$ of $\Bbb{Q}$, and we let $L_n$
denote the ring extension  of $\Bbb{Q}$; $L_n = \Bbb{Q}[t]/\langle c_n(t)\rangle$.
\end{definition}
We have 
$$L_n \cong  \bigoplus_{N\in\mathbb{N}_{\geq 3},N\mid n}\Bbb{Q}_N, \ n \geq 3.$$
(See Pilkington \cite{Me2}, Remarks 5.4 for details.)

The application of  Pilkington \cite{Me2}, Theorem 4.15  to the Coxeter systems associated to a  finite irreducible Coxeter group, gives the following  isomorphisms:

\begin{theorem} (Pilkington \cite{Me2}, Theorem 5.8) 
 Let $(W, S)$ be a Coxeter system with rank $N$, where $W$ is a finite irreducible Coxeter group.  We have the following unital ring isomomorphisms:
  $$\widehat{R}_{W,S} \cong \begin{cases} M_N(\Bbb{Q}) & \mbox{if $(W, S)$ is of type $A_n (n \geq 1), B_n (n \geq 2), D_n (n \geq 4), E_6, E_7, E_8, F_4$ or $G_2$,}\\
M_N(\Bbb{Q}_5) & \mbox{if $(W, S)$ is of type $H_3,$ or $H_4$,}\\
M_2(\Bbb{Q}_n) & \mbox{if $(W, S)$ is of type $I_2(n), n \geq 5, n \not= 6$.}\\
\end{cases}$$
We also have unital ring homomorphisms:
 $$\widehat{\widetilde{R}}_{W,S} \cong \begin{cases} M_N(\Bbb{Q}) & \mbox{if $(W, S)$ is of type $A_n (n \geq 1), D_n (n \geq 4), E_6, E_7$ or $E_8$,}\\
M_N(\Bbb{Q}_5) & \mbox{if $(W, S)$ is of type $H_3,$ or $H_4$,}\\
M_2({L}_n) & \mbox{if $(W, S)$ is of type  $I_2(n), n \geq 5$, where $n$ is odd.}\\
\end{cases}$$
\end{theorem}

For the remaining Weyl groups, with Coxeter systems of type $B_n \ (n \geq 2), F_4$ and $G_2 = I_2(6)$, 
and  the dihedral groups, with Coxeter systems of type   $I_2(n)$, where $n > 6$ and $n$ even, 
Pilkington \cite{Me2}, Theorem 4.15 does not apply directly to the algebras $\widehat{\widetilde{R}}_{W,S}$, due to the lack of an appropriate maximal tree in the graph $\Gamma_{W, S}$. 
 In these particular cases, we can show that  $\widehat{\widetilde{R}}_{W,S}$
 is isomorphic to 
 an algebra of type $R_{\mathfrak{F}_{Y_1}}$ associated to a graph $\Gamma = (S',Y')$, where $\Gamma$ is not connected and  $\mathfrak{F}_{Y_1}$ is a set of polynomials in $\Bbb{Q}[t]$
 with property F.  In Pilkington \cite{Me2}, Example 5.9,  a special case  of this isomorphism was demonstrated  for the  Coxeter system of type $B_4$. Here we 
 demonstrate the isomorphism for
 a general Coxeter System which encompasses all three of the exceptional types of Coxeter systems.

\subsection{ The Algebra  $\widehat{\widetilde{R}}_{W,S} $ associated to Coxeter Systems of Type $B_n \ (n \geq 2), F_4$,  $G_2$, 
and   $I_2(n)$, where $n > 6$.
}
Consider a Coxeter system,   $(W, S)$,  with $(n + m + 2) \times (n + m + 2)$ 
 Coxeter matrix of the form
   \[
 \left[\begin{array}{cccccccccccccc}
 1 &  3   & 2   & \cdots   &2 & 2 &2 & 2  & \cdots    &2 & 2& 2 \\
3 &   1       & 3 & \cdots  &2  & 2&2 & 2  & \cdots    &2 & 2& 2\\
2 & 3  & 1  & \cdots    & 2&2&2 & 2  & \cdots    &2 & 2& 2\\
\vdots  & \vdots  & \vdots     & \ddots     & \vdots&\vdots  & \vdots  & \vdots     & \ddots   & \vdots  & \vdots & \vdots  \\
2 & 2  & 2& \cdots    &1 & 3&2   &2& \cdots    &2 & 2& 2\\
2 & 2  &2& \cdots    &3 & 1 & a &2 & \cdots    &2 & 2& 2\\
2 & 2  &2& \cdots    &2 &a &   1       & 3& \cdots  &2  & 2& 2\\
2 & 2  &2& \cdots    &2 &2 & 3   & 1  & \cdots    & 2&2& 2\\
\vdots  & \vdots  & \vdots     & \ddots     & \vdots&\vdots  & \vdots  & \vdots     & \ddots   & \vdots  & \vdots & \vdots  \\
2 & 2  &2& \cdots    &2 & 2& 2  &2& \cdots    & 1&3&2\\
2 & 2  &2& \cdots    &2 & 2 & 2  & 2& \cdots   &3 &1 & 3\\
2 & 2  &2& \cdots    &2 & 2 & 2  &2& \cdots    &2&3& 1\\
\end{array}\right]
    \]

 where  $a$ is an even integer  greater than or equal to 4. 
 
 The associated graph, $\Gamma_{W, S} = (S, Y)$, where 
$S = \{s_1, s_2,  \dots , s_{n}, s_{n + 1}, t_{m+1},  t_{m}, \dots,  t_1\}$, with the order corresponding to the matrix entries,  is shown below.

\begin{equation*}
\xymatrix{
s_1\ar@{-}[r]^{3}&s_2\ar@{-}[r]^{3}&{s_3}\ar@{.}[r]&{s_{n}}\ar@{-}[r]^{3}&{s_{n + 1}}\ar@{-}[r]^{a}&{t_{m + 1}}\ar@{-}[r]^{3}&t_{m}\ar@{.}[r]&{t_{3}}\ar@{-}[r]^{3}&{t_2}\ar@{-}[r]^{3}&{t_1}\\
}
\end{equation*}
The geometric edges of the  graph $\Gamma_{W, S}$, $\left\{\{s_i, s_{i+1}\}, \{t_j, t_{j+1}\}, \{s_{n+1}, t_{m+1}\}, \  \vert \ 1 \leq i \leq n, 1 \leq j \leq m \right\},$  are labelled with 
the  corresponding entries from the Coxeter matrix. 
  The diagram  below represents the orientation  
  $$Y_+ =
  \{y_{i} = [s_i, s_{i + 1}], z = [s_{n + 1}, t_{m + 1}], x_j =  [t_j, t_{j + 1}] \ |\  1 \leq i \leq n, 1\leq j \leq m\} .$$
  
  \begin{equation*}
\xymatrix{
s_1\ar@{->}[r]^{y_{1}}&s_2\ar@{->}[r]^{y_{2}}&{s_3}\ar@{.}[r]&{s_{n}}\ar@{->}[r]^{y_{n}}&{s_{n+1}}\ar@{->}[r]^{z}&{t_{m+1}}\ar@{<-}[r]^{x_{m}}&t_{m}\ar@{.}[r]&{t_{3}}\ar@{<-}[r]^{x_{2}}&{t_2}\ar@{<-}[r]^{x_1}&{t_1}\\
}
\end{equation*}
We have
  $Y = Y_+ \bigcupdot \overline{Y_+}$.

We have that $\widehat{\widetilde{I}}_{W,S}$ is the two-sided ideal of $\widehat{P}_{W,S}$ generated by the set of elements
\begin{equation*}
 \left\{
 \begin{aligned}
 c_{3}(y_{i}\overline{y_{i}}) = 
y_{i}\overline{y_{i}} - [s_i], \\
c_{3}(\overline{y_{i}}y_{i}) = 
\overline{y_{i}}y_{i} - [s_{i + 1}],
\end{aligned} \ \ \ \ 
 \begin{aligned}
\overline{z}c_{a}(z\overline{z}) = \overline{z}(h(z\overline{z}) - k) =  \overline{z}(z\overline{z}h_a(z\overline{z}) - k[s_{n+1}]),\\
zc_{a}(\overline{z}z) = z(h(\overline{z}z) - k[t_{m+1}]) = z(\overline{z}zh_a(\overline{z}z) - k),
\end{aligned} \ \ \ \ 
 \begin{aligned}
c_{3}(x_{j}\overline{x_{j}}) 
= x_{j}\overline{x_{j}} - [s_j],\\
c_{3}(\overline{x_{j}}x_{j}) 
= \overline{x_{j}}x_{j} - [s_{j + 1}]
\end{aligned}\right\}
\end{equation*} 
for  $1 \leq i \leq n, 1 \leq j \leq m$.

Let $\Gamma = (S',Y')$ be the graph,  shown in Figure \ref{fig:fig2}, below with vertices 
$$S' = \left\{s'_i,  t'_j | 1\leq i \leq n+1, 1\leq j \leq m+1\right\} \ \bigcup
\ \left\{s''_i,  t''_j | 1\leq i \leq n+1, 1\leq j \leq m+1\right\}$$
 and edges $Y' = Y'_+ \bigcupdot \overline{Y'_+}$ where 
$ Y'_+  = \left\{ y'_i , y''_i, x'_j, x''_j\ | 1 \leq i \leq n,  1 \leq j \leq m\right\}  \bigcup \{z'\}$(see Figure \ref{fig:fig3} below).  

\begin{figure}[H]
$$
\xymatrix{
s'_1\ar@{-}[r]^{3}&s'_2\ar@{-}[r]^{3}&{s'_3}\ar@{.}[r]&{s_{n}}\ar@{-}[r]^{3}&{s'_{n + 1}}\ar@{-}[r]^{a}&{t'_{m + 1}}\ar@{-}[r]^{3}&t'_{m}\ar@{.}[r]&{t'_{3}}\ar@{-}[r]^{3}&{t'_2}\ar@{-}[r]^{3}&{t'_1}\\\
s''_1\ar@{-}[r]^{3}&s''_2\ar@{-}[r]^{3}&{s''_3}\ar@{.}[r]&{s_{n}}\ar@{-}[r]^{3}&{s''_{n + 1}}&{t''_{m + 1}}\ar@{-}[r]^{3}&t''_{m}\ar@{.}[r]&{t''_{3}}\ar@{-}[r]^{3}&{t''_2}\ar@{-}[r]^{3}&{t'_1}\\
}
$$
\caption{Graph of $\Gamma = (S',Y')$.}
\label{fig:fig2}
\end{figure}
  The edges corresponding to the orientation  $Y'_+$ are shown in Figure \ref{fig:fig3} below.
For $1 \leq i \leq n$, we have  $o(y'_i) = s'_i$,  $t(y'_i) = s'_{i + 1}$,  $o(y''_i) = s''_i$ and $t(y''_i) = s''_{i + 1}$. 
For $1 \leq j \leq m$, we have $o(x'_j) = t'_j$, $t(x'_j) = t'_{j+1}$, 
$o(x''_j) = t''_j$ and $t(x''_j) = t''_{j+1}$.
We have  $o(z') = s'_{n+1}$,  $t(z') = t'_{n+1}$.
\begin{figure}[H]
$$
\xymatrix{
s'_1\ar@{->}[r]^{y'_{1}}&s'_2\ar@{->}[r]^{y'_{2}}&{s'_3}\ar@{.}[r]&{s_{n}}\ar@{->}[r]^{y'_{n}}&{s'_{n + 1}}\ar@{->}[r]^{z'}&{t'_{m + 1}}\ar@{<-}[r]^{x'_{m}}&t'_{m}\ar@{.}[r]&{t'_{3}}\ar@{<-}[r]^{x'_{2}}&{t'_2}\ar@{<-}[r]^{x'_{1}}&{t'_1}\\\
s''_1\ar@{->}[r]^{y''_{1}}&s''_2\ar@{->}[r]^{y''_{2}}&{s''_3}\ar@{.}[r]&{s_{n}}\ar@{->}[r]^{y''_{n}}&{s''_{n + 1}}&{t''_{m + 1}}\ar@{<-}[r]^{x''_{m}}&t''_{m}\ar@{.}[r]&{t''_{3}}\ar@{<-}[r]^{x''_{2}}&{t''_2}\ar@{<-}[r]^{x''_{1}}&{t'_1}\\
}
$$
\caption{The orientation $Y'_+$ on the graph  $\Gamma = (S',Y')$.}
\label{fig:fig3}
\end{figure}

  The graph $\Gamma = (S', Y')$,  has three connected components, $\Gamma_{\dot{S}, \dot{Y}}, \Gamma_{\ddot{S}, \ddot{Y}}$ and $\Gamma_{\dddot{S}, \dddot{Y}}$,
where the vertices are given by 
$\dot{S} = \{s_i', t'_j \ | \  1\leq i \leq n+1, 1 \leq j \leq m + 1\}$, 
$\ddot{S} = \{s''_i \ | \ 1 \leq i \leq n+1\}$ and 
$\dddot{S} =  \{t''_j \ | \ 1 \leq j \leq m+1\}$ respectively and 
the edges are given by 
 $\dot{Y} = \{y_i', \overline{y'_i},  x_j', \overline{x'_j}, z', \overline{z'}\ |\  1\leq i \leq n, 1\leq j \leq m\}$,
 $ \ddot{Y} = \{y_i'',  \overline{y_i''}\ |\ 1 \leq i \leq n\}$ and $\dddot{Y} = \{x_j'',  \overline{x_j''}\ |\  1\leq j \leq m\}$ respectively.

  We let $Y_1 = Y'$ and we let $ \frak{F}_{Y_1} = \left\{f_y(t) \ |\ y \in Y_1 \right\}$, where $f_y(t) = f_{\overline{y}}(t)$,  is the set of polynomials such that
  \begin{equation*}
  \begin{aligned}
 f_{y'_i}(t)  &= f_{x'_j}(t)  &=& \ C_3(t), &\ \ \mbox{for} & \ 1\leq i \leq n, 1 \leq j \leq m,\ \\
  f_{y''_i}(t) &= f_{x''_j}(t)  &=& \ C_3(t), &\ \ \mbox{for} & \ 1\leq i \leq n, 1 \leq j \leq m,\ \\
f_{z'}(t)   &= c_a(t) &=& \prod_{{\tiny \left.\begin{array}{c}N\in\mathbb{N}_{\geq 3} \\N\mid a\end{array}\right.}}C_N(t).&\\
\end{aligned}
\end{equation*}
The set of polynomials  $\frak{F}_{Y_1}$ has property F by Remark \ref{seen}.  By
Remarks \ref{relations} and 
Theorem 
\ref{theisocor}, we  get a unital ring isomorphism 
$$R_{\frak{F}_{Y_1}}  
 \cong M_{m + n + 2}({L}_a) \bigoplus M_{n}(\Bbb{Q}) \bigoplus M_{m}(\Bbb{Q}).$$

 \begin{theorem}{\label{general}}
Let the graphs $\Gamma_{W, S}$ and  $\Gamma$, and their corresponding algebras $\widehat{\widetilde{R}}_{W,S}$  and $R_{\frak{F}_{Y_1}} $ be as defined above. 
Then we have unital ring isomorphisms: 
$$\widehat{\widetilde{R}}_{W,S} \cong R_{\frak{F}_{Y_1}} \cong 
  M_{m + n + 2}({L}_a) \bigoplus M_{n}(\Bbb{Q}) \bigoplus M_{m}(\Bbb{Q}).$$
 \end{theorem}
\begin{proof}
We let $p^{\star}$ denote the coset  $p + \widehat{\widetilde{I}}_{W, S} \in \widehat{\widetilde{R}}_{W, S}$ for $p \in \widehat{P}_{W, S}$, and we let 
$q^{\dagger}$
denote the coset $q + I_{\mathfrak{F}_{Y_1}} \in R_{\mathfrak{F}_{Y_1}}$ for $q \in 
\Bbb{Q}\Gamma$ throughout. 
We define  homomorphisms $\Theta: \widehat{\widetilde{R}}_{W, S} \to R_{\mathfrak{F}_{Y_1}}$
and $\Omega: R_{\mathfrak{F}_{Y_1}} \to \widehat{\widetilde{R}}_{W, S}$ and show that 
$\Theta = \Omega^{-1}$, which will prove the result.

 \noindent
{\bf The Homomorphism \underline{$\displaystyle \Theta: \widehat{\widetilde{R}}_{W, S} \to R_{\mathfrak{F}_{Y_1}}$}.} 
We let $\theta :\widehat{P}_{W, S} \to R_{\mathfrak{F}_{Y_1}}$ be defined on the vertices and edges of $\Gamma_{W, S}$ as follows:

\begin{equation*}
\begin{aligned}
\theta([s_i]) &= [s'_i]^\dagger + [s''_i]^\dagger& \mbox{ for} \  1 \leq i \leq n+1, \\
\theta([t_j]) &= [t'_j]^\dagger + [t''_j]^\dagger& \mbox{ for} \  1 \leq i \leq m+1, \\
\theta(y_i) &= (y'_i)^\dagger + (y''_i)^\dagger, & \mbox{ for } \   1 \leq i \leq n, \\
\theta(x_j) &= (x'_j)^\dagger + (x''_j)^\dagger, & \mbox{ for } \   1 \leq j \leq m, \\
\end{aligned}\\
\hskip 1in
\begin{aligned}
\theta(z) &= (z')^\dagger,&\\
\theta(\overline{z}) &=  (\overline{z'})^\dagger,&\\
\theta(\overline{y_i}) &= (\overline{y'_i})^\dagger + (\overline{y''_i})^\dagger, & \mbox{ for } \   1 \leq i \leq n, \\
\theta(\overline{x_j}) &= (\overline{x'_j})^\dagger + (\overline{x''_j})^\dagger, & \mbox{ for } \   1 \leq j \leq m, \\
\end{aligned}\\
\end{equation*}
It is not difficult to see that the set $\{\theta(s_i), \theta(t_j)\ | \ 1 \leq i \leq n + 1, \ 1 \leq j \leq m+1\}$ is a set of  mutually orthogonal  idempotents in
$R_{\mathfrak{F}_{Y_1}}$. It is also easy to check that products of vertices and edges  in 
 $\widehat{P}_{W, S}$ 
are mapped to the product of their images by  $\theta$. Therefore, $\theta$ is a well defined homomorphism 
from $\widehat{P}_{W, S}$  to  $R_{\mathfrak{F}_{Y_1}}$. 
Applying $\theta$ to the generators of $\widehat{\widetilde{I}}_{W, S}$. 
We see that for $1 \leq i \leq n$; 
\begin{equation*}
\begin{aligned}
\theta(y_i\overline{y_i} - [s_i]) 
&= \left((y'_i)^\dagger + (y''_i)^\dagger\right) \left((\overline{y'_i})^\dagger + (\overline{y''_i})^\dagger\right)
- ([s'_i]^\dagger +[s''_i]^\dagger)  \\
&= (y'_i)^\dagger(\overline{y'_i})^\dagger - [s'_i]^\dagger + (y''_i)^\dagger(\overline{y''_i})^\dagger - [s''_i]^\dagger\\
&= 0. 
\end{aligned}
\end{equation*}
Similar calculations give that   $\theta(\overline{y_i} y_i - [s_{i + 1}])  = 0$ for $1 \leq i \leq n$,   
$\theta(x_j\overline{x_j} - [t_{j}])  = 0$ and 
 $\theta(\overline{x_j}x_j - [t_{j+1}])  = 0$ for $1 \leq i \leq m$.
 
 \noindent
We have 
\begin{equation*}
\begin{aligned}
\theta\left(\overline{z}c_a(z\overline{z})\right) 
&= (\overline{z'})^\dagger\left(c_a(z'\overline{z'})\right)^\dagger
&= 0,\\
\theta\left(zc_a(\overline{z}z)\right) 
&= ({z'})^\dagger\left(c_a(\overline{z'}z')\right)^\dagger
&= 0.
\end{aligned}
\end{equation*}
Therefore, we have $\widehat{\widetilde{I}}_{W, S} \subseteq \ker \theta$, 
and the quotient homomorphism $\Theta : \widehat{\widetilde{R}}_{W, S} \to R_{\frak{F}_{Y_1}}$
which sends $p^\star$ to $\theta(p)$ for $p \in \widehat{P}_{W, S}$, is well defined.

\noindent
{\bf The Homomorphism \underline{$\displaystyle \Omega: R_{\mathfrak{F}_{Y_1}} \to  \widehat{\widetilde{R}}_{W, S}$}.}
Let  $\omega: \Bbb{Q}\Gamma \to \widehat{\widetilde{R}}_{W,S}$ be defined on the vertices and edges of $\Gamma$ as follows:
\begin{equation*}
\begin{aligned}
\omega([s'_{n+1}]) &= h(z\overline{z})^{\star}/{k},\\
\omega([t'_{m+1}]) &= h(\overline{z}z)^{\star}/{k},\\
\end{aligned}
\hskip .5in
\begin{aligned}
 \omega([s''_{n + 1}]) &= [s_{n + 1}]^{\star} - \omega([s'_{n + 1}]),\\
  \omega([t''_{m + 1}]) &= [t_{m + 1}]^{\star} - \omega([t'_{m + 1}]),
\end{aligned}
\hskip .5in
\begin{aligned}
 \omega(z') &= z^\star,\\
   \omega(\overline{z'}) &= \overline{z}^\star,
\end{aligned}
\end{equation*}
and for $1 \leq i \leq n$ and $1 \leq j \leq m$, we let 
\begin{equation*}
\begin{aligned}
\omega([s_i']) &=  y_i^\star\omega([s'_{i+1}])\overline{y_i}^\star, \\
\omega([t_j']) &=  x_j^\star\omega([t'_{j+1}])\overline{x_j}^\star, \\
 \omega(y'_i) &= \omega([s'_i])y_i^\star\omega([s'_{i+1}]),\\
  \omega(y''_i) &= \omega([s''_i])y_i^\star\omega([s''_{i+1}])\\
   \omega(x'_j) &= \omega([t'_j])x_j^\star\omega([t'_{j+1}]),\\
  \omega(x''_j) &= \omega([t''_j])x_j^\star\omega([t''_{j+1}]),
\end{aligned}
\hskip .5in
\begin{aligned}
\omega([s''_i]) &=  y_i^\star\omega([s''_{i+1}])\overline{y_i}^\star, \\
\omega([t''_j]) &=  x_j^\star\omega([t''_{j+1}])\overline{x_j}^\star, \\
 \omega(\overline{y'_i}) &= \omega([s'_{i+1}])\overline{y_i}^\star\omega([s'_i]),\\
 \omega(\overline{y''_i}) &= \omega([s''_{i+1}])\overline{y_i}^\star\omega([s''_i]),\\
  \omega(\overline{x'_j}) &= \omega([t'_{j+1}])\overline{x_j}^\star\omega([t'_j]),\\
 \omega(\overline{x''_j}) &= \omega([t''_{j+1}])\overline{x_j}^\star\omega([t''_j]).
\end{aligned}
\end{equation*}

In order to show that $\omega$ gives  a well defined homomorphism from $\Bbb{Q}\Gamma$ to $\widehat{\widetilde{R}}_{W,S}$, we need to show that the map preserves the relations in the path algebra. 
In our calculations, we will use the  relations on the images of the edges in
$\widehat{\widetilde{R}}_{W, S}$ listed below;
\begin{equation}\label{rellies2}
 \begin{aligned}
 \left(c_{3}(y_{i}\overline{y_{i}})\right)^\star = 
y_{i}^\star\overline{y_{i}}^\star - [s_i]^\star = 0, \\
\overline{z}^\star\left(c_{a}(z\overline{z})\right)^\star = \overline{z}^\star(h(z\overline{z}) - k[s_{n+1}])^\star = 0,\\
\left(c_{3}(x_{j}\overline{x_{j}})\right)^\star
= x_{j}^\star\overline{x_{j}}^\star - [s_j]^\star = 0,\\
\end{aligned} \ \ \ \ 
 \begin{aligned}
 \left(c_{3}(\overline{y_{i}}y_{i})\right)^\star = 
\overline{y_{i}}^\star y_{i}^\star - [s_{i + 1}]^\star = 0,\\
z^\star\left(c_{a}(\overline{z}z)\right)^\star = z^\star(h(\overline{z}z) - k[t_{m+1}])^\star = 0,\\
\left(c_{3}(\overline{x_{j}}x_{j})\right)^\star
= \overline{x_{j}}^\star x_{j}^\star - [s_{j + 1}]^\star = 0,
\end{aligned} \ \ \ \ 
\end{equation}
where $1 \leq i \leq n$ and $1 \leq j \leq m$. \\
The following remarks will also be useful throughout our calculations:
\begin{Remark} \label{sis2}

 \begin{enumerate}
\item For  $1 \leq i \leq n+1$, $1 \leq j \leq m+1$,  
\begin{equation*}
\begin{aligned}
[s_i]^\star\omega([s_i']) &= \omega([s_i'])[s_i]^\star &= \omega([s_i']),\\
[s_i]^\star\omega([s''_i]) &= \omega([s''_i])[s_i]^\star &= \omega([s''_i]),
\end{aligned}
\hskip .75in
\begin{aligned}
[t_j]^\star\omega([t_j']) &= \omega([t_j'])[t_j]^\star &= \omega([t_j']),\\
[t_j]^\star\omega([t''_j]) &= \omega([t''_j])[t_j]^\star &= \omega([t''_j]).
\end{aligned}
\end{equation*}
This follows 
from the definition of $\omega([s'_i])$,  $\omega([s''_i])$, $\omega([t'_j])$ and $\omega([t''_j])$\\
\item For  
$1 \leq i \leq n + 1$  we have that 
\begin{equation}\label{siprimeprime}
\omega([s''_i]) = [s_i]^\star - \omega([s'_i])
\end{equation} 
This is true by definition for $i = n + 1$. Using induction, we assume that Equation \ref{siprimeprime} holds  for 
all $l$ such that $i + 1 \leq l \leq n+1$. 
We then have
 \begin{equation*}
 \begin{aligned}
\omega([s''_i]) &= y_i^\star\omega([s''_{i+1}])\overline{y_i}^\star\\
 &=y_i^\star([s_{i+1}]^\star - \omega([s'_{i+1}]))\overline{y_i}^\star\\
 &= y_i^\star[s_{i+1}]^\star\overline{y_i}^\star - y_i^\star\omega([s'_{i+1}])\overline{y_i}^\star\\
 &= [s_i]^\star - \omega([s'_i]) \\
 & \mbox{(from the definition of $\omega([s'_i])$ and the relations listed in (\ref{rellies2}).}
 \end{aligned}
 \end{equation*}
 A similar calculation gives that  
 $$\omega([t''_j]) = [t_j]^\star - \omega([t'_j])$$
  for 
$1 \leq j \leq m + 1$.
\\
\item 
$$\left(\left(h(z\overline{z})\right)^{\star}\right)^2 = k\left(h(z\overline{z})\right)^{\star}.$$
We have 
$c_a(t) = h(t) - k = th_a(t) - \kappa$ where $h_a(t) \in \Bbb{Z}[t]$. From \ref{rellies2} we have  identities
$0 = \left(\overline{z}c_{a}(z\overline{z})\right)^{\star} = \left(\overline{z}(h(z\overline{z}) - k[s_{n+1}])\right)^{\star}$ and 
$0 = \left(zc_{a}(\overline{z}z)\right)^{\star} = \left(z(h(\overline{z}z) - k[t_{m+1}])\right)^{\star}$  in the algebra 
$\widehat{\widetilde{R}}_{W,S}$.  Thus we have the following:
\begin{equation*}
 \left(\left(h(z\overline{z})\right)^{\star}\right)^2 = \left(z\overline{z}h_a(z\overline{z})h(z\overline{z})\right)^{\star} = \left(h_a(z\overline{z})z\overline{z}h(z\overline{z})\right)^{\star} = \left(h_a(z\overline{z})z\overline{z}k[s_{n+1}]\right)^{\star} = k\left(h(z\overline{z})\right)^{\star}.
\end{equation*}
Likewise we have that 
\begin{equation*}
\left(\left(h(\overline{z}z)\right)^{\star}\right)^2  = k\left(h(\overline{z}z)\right)^{\star}.
\end{equation*}
 \end{enumerate}
 \end{Remark}

We first show  that $\{\omega([s_i']), \omega([s_i'']), \omega([t_j']), \omega([t_j'']) \ |\  1 \leq i \leq n + 1, 1 \leq j \leq m+1\}$, are a set of 
mutually orthogonal idempotents in
$\widehat{\widetilde{R}}_{W, S}$. Using the definition of $\omega$ and Remarks \ref{sis2}(3),
we see that 
$$(\omega([s'_{n + 1}]))^2 = \left(\left(h(z\overline{z})\right)^{\star}\right)^2/{k^2}
=k\left(h(z\overline{z})\right)^{\star}/{k^2} = \left(h(z\overline{z})\right)^{\star}/{k} =  \omega([s'_{n+1}]).$$
By a similar calculation, we have that $\omega([t'_{m+1}])$ is also an idempotent.

We use induction to show that $\omega([s'_i])$ is an idempotent for $1 \leq i \leq n$. 
Assuming that $\omega([s'_l])$ is an idempotent for all $l$ such that $i + 1 \leq l \leq n + 1$, then we have 
\begin{equation*}
\begin{aligned}
(\omega([s'_i]))^2 
&= y_i^\star\omega([s'_{i+1}])\overline{y_i}^\star y_i^\star\omega([s'_{i+1}])\overline{y_i}^\star
 \ \mbox{(from the definition of $\omega$)}\\
 &= y_i^\star\omega([s'_{i+1}])[s_{i+1}]^\star\omega([s'_{i+1}])\overline{y_i}^\star
  \ \mbox{(using the relations listed in  (\ref{rellies2}) )}\\
&= y_i^\star\omega([s'_{i+1}])\overline{y_i}^\star 
  \ \mbox{(from Remarks \ref{sis2} (1) )}\\
  &= \omega([s'_i]) 
 \ \mbox{(from the definition of $\omega$)}.\\
\end{aligned}
\end{equation*}
In a similar way, we see that $\omega([t'_j])$ is an idempotent for $1 \leq j \leq m$.

Applying  Remarks \ref{sis2} (1) and (2) to $\omega(s''_i)$, where $1 \leq i \leq n +1$, we get 
$$(\omega(s''_i))^2 = ([s_i]^\star - \omega(s'_i))^2 = [s_i]^\star -[s_i]^\star\omega(s'_i) - \omega(s'_i) [s_i]^\star + (\omega(s'_i))^2 = [s_i]^\star - \omega(s'_i) = \omega(s''_i).$$
Similarly, $(\omega(t''_j))^2 = \omega(t''_j)$ for $1 \leq j \leq m +1$.
Using Remarks \ref{sis2} (2), we show  that for  a given $i$ with $1 \leq i \leq n+1$, 
$\omega([s'_i])$ and $\omega([s''_i])$ are orthogonal;
\begin{equation*}
\begin{aligned}
\omega([s'_i])\omega([s''_i]) 
&= \omega([s'_i])([s'_i] - \omega([s'_i]))\\
&= \omega([s'_i]) - (\omega([s'_i]))^2 \ \mbox{(by Remarks \ref{sis2} (1))}\\
&= \omega([s'_i]) - \omega([s'_i]) = 0.\\
\end{aligned}
\end{equation*}
Similarly we get that $\omega([s''_i])\omega([s'_i]) = 0$ and therefore, $\omega([s'_i])$ and $\omega([s''_i])$ are orthogonal idempotents for any given $i$ with $1 \leq i \leq n+1$.
By symmetry or by similar routine calculations, we see that $\omega([t'_j])$ and $\omega([t''_j])$ 
are orthogonal idempotents for $1 \leq j \leq m+1$. 
The remaining verifications of mutual orthogonality for the set of idempotents, 
$$\frak{S} = \{\omega([s_i']), \omega([s_i'']), \omega([t_j']), \omega([t_j'']) \ |\  1 \leq i \leq n + 1, 1 \leq j \leq m+1\},$$  can be checked by comparing the origin and terminus of the paths in the coset representatives 
used in  their definitions.

That the product of a given vertex  and edge in $\Bbb{Q}\Gamma$ is mapped to the product of the images in $\widehat{\widetilde{R}}_{W,S}$ by $\omega$ now follows easily from the definition of $\omega$ on the edges, and the fact that the set  $\frak{S}$ above is a set of mutually orthogonal idempotents. 
Therefore, we see that $\omega$ gives us a 
well defined homomorphism  $\omega: \Bbb{Q}\Gamma \to \widehat{\widetilde{R}}_{W,S}$.

Looking at  at the behavior of $\omega$ on the generators of $I_{\frak{F}_{Y_1}}$, 
we have 
$$\omega(c_a(z'\overline{z'})) = \omega(h(z'\overline{z'}) - k[s'_{n+1}])
= \left(h(z\overline{z})\right)^\star -  k\left(h(z\overline{z})\right)^\star/{k} = 0.$$
Similarly $\omega(c_a(\overline{z'}z')) = 0$.

For $1 \leq i \leq n$, we have 
\begin{equation*}
\begin{aligned}
\omega(C_3(y_i'\overline{y'_i})) &= \omega((y'_i\overline{y'_i}) - [s'_i])\\
&= \omega([s'_i])y^\star_i\omega([s'_{i+1}])\omega([s'_{i+1}])\overline{y_i}^\star\omega([s'_i])
 -  
\omega([s'_i])\\
&= \omega([s'_i])y^\star_1\omega([s'_{i+1}])\overline{y_i}^\star\omega([s'_i])
 -  
\omega([s'_i])\\
&= \omega([s'_i])\omega([s'_i])\omega([s'_i]) 
 -  
\omega([s'_i]) \ \mbox{(from the definition of $\omega([s'_i])$)}\\
 &= 0\\
\end{aligned}
\end{equation*}
Similar calculations yield that $\omega(C_3(y''_i\overline{y''_i})) = 0$, 
$\omega(C_3(x'_j\overline{x'_j})) = 0$ and  $\omega(C_3(x''_j\overline{x''_j})) = 0$ 
for $1 \leq i \leq n, 1 \leq j \leq m$.

For $1 \leq i \leq n$, we have 
\begin{equation*}
\begin{aligned}
\omega(C_3(\overline{y'_i} y_i')) &= \omega((\overline{y'_i} y'_i) - [s'_{i+1}])\\
&= \omega([s'_{i+1}])\overline{y_i}^\star\omega([s'_i])\omega([s'_i]){y}^\star_i\omega([s'_{i+1}])
 -  \omega([s'_{i+1}])\\
&= \omega([s'_{i+1}])\overline{y_i}^\star\omega([s'_i]){y}^\star_i\omega([s'_{i+1}])
 -  \omega([s'_{i+1}])\\
&= \omega([s'_{i+1}])\overline{y_i}^\star{y}^\star_i\omega([s'_{i+1}])\overline{y_i}^\star{y}^\star_i\omega([s'_{i+1}])
 -  \omega(s'_{i+1})\ \mbox{(from the definition of $\omega([s'_i])$)}\\
&=\omega([s'_{i+1}])^3  -  \omega([s'_{i+1}]) \ \mbox{(from the relations listed in  (\ref{rellies2}))}\\
& =0 \\
\end{aligned}
\end{equation*}
Similar  calculations yield that $\omega(C_3(\overline{y''_i} y''_i)) = 0$,
$\omega(C_3(\overline{x'_j}x'_j)) = 0$ and  $\omega(C_3(\overline{x''_j}x''_j)) = 0$ 
for $1 \leq i \leq n, 1 \leq j \leq m$.

This gives us that $I_{\frak{F}_{Y_1}}  \subseteq \ker \omega$ and that the quotient map 
$\Omega: R_{\frak{F}_{Y_1}} \to  \widehat{\widetilde{R}}_{W, S}$ which sends $p + I_{\frak{F}_{Y_1}}$
to $\omega(p)$ in $ \widehat{\widetilde{R}}_{W, S}$ is well defined.

\noindent
{\bf  The  Homomorphism \underline{$\displaystyle \Omega \circ \Theta: \widehat{\widetilde{R}}_{W, S} \to \widehat{\widetilde{R}}_{W, S}$ }. }
We check that the map $\Omega \circ \Theta$ is the identity map on on the generators of  $\widehat{\widetilde{R}}_{W, S}$.
We have 
$\Omega \circ \Theta(z^\star) = \Omega\left((z')^\dagger\right) = z^\star.$ Similarly,
$\Omega \circ \Theta(\overline{z}^\star) = \Omega\left((\overline{z'})^\dagger\right) = \overline{z}^\star.$

For $1 \leq i \leq n$, we have 
\begin{equation*}
\begin{aligned}
\Omega \circ \Theta(y_i^\star) 
&= \Omega\left((y_i')^\dagger\right) + \Omega\left((y_i'')^\dagger\right) \\
&= \omega([s_i'])y_i^\star\omega([s_{i+1}']) +  \omega([s_i''])y_i^\star\omega([s_{i+1}'']) \\
&= y_i^\star\omega([s_{i+1}'])\overline{y_i}^\star y_i^\star \omega([s_{i+1}']) 
+ y_i^\star\omega([s_{i+1}''])\overline{y_i}^\star y_i^\star \omega([s_{i+1}'']) \ \mbox{(from the definition of  $\omega([s_i'])$ and $ \omega([s''_i])$)}\\
&= y_i^\star\omega([s_{i+1}'])[s_{i+1}]^\star\omega([s_{i+1}']) + y_i^\star\omega([s_{i+1}''])[s_{i+1}]^\star\omega([s_{i+1}'']) \ \ \mbox{(by the relations listed in (\ref{rellies2}))}\\
&= y_i^\star\omega([s_{i+1}']) + y_i^\star\omega([s_{i+1}''])  \ \mbox{(by Remarks \ref{sis2} (1))}\\
&= y_i^\star(\omega([s_{i+1}']) + [s_{i+1}]^\star - \omega([s_{i+1}'])) \ \mbox{(by Remarks \ref{sis2} (2))} \\
&= y_i^\star[s_{i+1}]^\star\\
  &= y_i^\star.
\end{aligned}
\end{equation*}
Similarly, using the relations (\ref{rellies2})) and Remarks \ref{sis2}, we get for $1 \leq i \leq n$
\begin{equation*}
\begin{aligned}
\Omega \circ \Theta(\overline{y_i}^\star) 
&= \Omega\left((\overline{y_i'})^\dagger\right) + \Omega\left((\overline{y_i''})^\dagger\right) \\
&= \omega([s_{i+1}'])\overline{y_i}^\star\omega([s_i']) +  \omega([s_{i+1}''])\overline{y_i}^\star\omega([s_i'']) \\
&= \omega([s_{i+1}'])\overline{y_i}^\star y_i^\star\omega([s_{i+1}'])\overline{y_i}^\star 
 +  \omega([s_{i+1}''])\overline{y_i}^\star y_i^\star\omega([s_{i+1}''])\overline{y_i}^\star  \\
&=  \omega([s_{i+1}'])[s_{i+1}]^\star\omega([s_{i+1}'])\overline{y_i}^\star 
 +  \omega([s_{i+1}''])[s_{i+1}]^\star\omega([s_{i+1}''])\overline{y_i}^\star  \\
&= \omega([s_{i+1}'])\overline{y_i}^\star + \omega([s_{i+1}''])\overline{y_i}^\star  \\
&= [s_{i+1}]^\star\overline{y_i}^\star  = \overline{y_i}^\star.  \\\end{aligned}
\end{equation*}

\noindent
By symmetry (or by a similar calculation), one gets that $\Omega \circ \Theta(x_j^\star) 
= x_j^\star$ and $\Omega \circ \Theta(\overline{x_j}^\star) = \overline{x_j}^\star$ for $1 \leq j \leq m$. 

\noindent
Using Remarks \ref{sis2} (2), we see  that for  $1 \leq i \leq n+1$, we have
$$\Omega \circ \Theta([s_i]^\star) = \Omega([s'_i]^\dagger + [s''_i]^\dagger)
= \omega([s'_i]) +  \omega([s''_i]) = \omega([s'_i])  + [s_i]^\star - \omega([s'_i])  = [s_i]^\star.$$
Similarly $\Omega \circ \Theta([t_j]^\star) =  [t_j]^\star$ for $1 \leq j \leq m$.

\noindent
Thus, since $\Omega \circ \Theta$ is a homomorphism, which acts as the identity on all of the generators of $\widehat{\widetilde{R}}_{W, S}$, it is the identity homomorphism on that 
algebra.

{\bf  The Homomorphism \underline{$\displaystyle \Theta\circ \Omega : R_{\frak{F}_{Y_1}} \to R_{\frak{F}_{Y_1}}$}. }
We check that $\Theta\circ \Omega$ acts as the identity on the generators of $R_{\frak{F}_{Y_1}}$.

First we look at the action of $\Theta\circ \Omega$ on the cosets associated to the vertices of the graph $\Gamma$, 
$$S' = \displaystyle \left\{[s'_i]^\dagger,  [s''_i]^\dagger, [t'_j]^\dagger,  [t''_j]^\dagger \ |\ 1\leq i \leq n+1, 1 \leq j \leq m+1 \right\}.$$
We have,
$$\Theta\circ \Omega([s'_{n + 1}]^\dagger) = \Theta\left(\left(h(z\overline{z})\right)^\star/k\right) = \left(h(z'\overline{z'})\right)^\dagger /k = k[s'_{n + 1}]^\dagger/k = [s'_{n + 1}]^\dagger.$$
By a similar calculation, we have $\Theta\circ \Omega([t'_{m+1}]^\dagger) =  [t'_{m+1}]^\dagger.$

\noindent
For $s'_i, 1 \leq i \leq n$ we use induction. Assuming that  $\Theta\circ \Omega([s'_{l}]^\dagger) =  [s'_{l}]^\dagger$ for $i + 1 \leq l \leq n+1$, we 
have that
$$\Theta\circ \Omega([s'_i]^\dagger) = \Theta(y_i^\star\omega([s'_{i+1}])\overline{y_i}^\star)
= (y'_i + y''_i)^\dagger[s'_{i+1}]^\dagger(\overline{y'_i} + \overline{y''_i})^\dagger = 
(y'_i)^\dagger[s'_{i+1}]^\dagger(\overline{y'_i})^\dagger = [s'_i]^\dagger.$$
A similar calculation gives that $\Theta\circ \Omega([t'_j]^\dagger) = [t'_j]^\dagger$ for 
$1 \leq j \leq m$. 

\noindent
 For $1 \leq i \leq n+1$, by Remarks \ref{sis2} (2), we have
$$\Theta\circ \Omega([s''_i]^\dagger) = \Theta([s_i]^\star - \omega([s'_i]^\star)) = [s'_i]^\dagger + [s''_i]^\dagger - \Theta\circ \Omega([s'_i]^\dagger) = [s'_i]^\dagger + [s''_i]^\dagger - [s'_i]^\dagger
= [s''_i]^\dagger.$$
Similarly, we see that $\Theta\circ \Omega([t''_j]^\dagger) =  [t''_j]^\dagger$ for $1 \leq j \leq m+1$.

Next we look at the  action of $\Theta\circ \Omega$ on the cosets associated to the edges of
 $\Gamma$.
We get that
$$\Theta\circ \Omega((z')^\dagger) = \Theta(z^\star) = (z')^\dagger  \ \mbox{and}  \ \Theta\circ \Omega((\overline{z'})^\dagger) = \Theta(\overline{z}^\star) = (\overline{z'})^\dagger.$$
For $1 \leq i \leq n$, we have 
\begin{equation*}
\begin{aligned}
\Theta\circ \Omega((y'_i)^\dagger) 
&= \Theta\left(\omega([s'_i])y_i^\star\omega([s'_{i+1}])\right)\\
&= \left(\Theta\circ \Omega([s'_i]^\dagger)\right)((y'_i)^\dagger + (y''_i)^\dagger)\left(\Theta\circ \Omega([s'_{i+1}]^\dagger)\right)\\
& = [s'_i]^\dagger((y'_i)^\dagger + (y''_i)^\dagger)[s'_{i+1}]^\dagger\\
&=  [s'_i]^\dagger(y'_i)^\dagger[s'_{i+1}]^\dagger = (y'_i)^\dagger. 
\end{aligned}
\end{equation*}
Similar calculations show that for $1 \leq i \leq n$,  $\Theta\circ \Omega((\overline{y'_i})^\dagger)  = (\overline{y'_i})^\dagger$, 
$\Theta\circ \Omega(({y}''_i)^\dagger)  = ({y}''_i)^\dagger$ and $\Theta\circ \Omega((\overline{y''_i})^\dagger)  = (\overline{y''_i})^\dagger$.
Using symmetry, or performing similar calculations, for $1 \leq j \leq m$, we have  $
\Theta\circ \Omega(({x}'_j)^\dagger)  = ({x}'_j)^\dagger$, 
$\Theta\circ \Omega((\overline{x'_j})^\dagger)  = (\overline{x'_j})^\dagger$, 
$\Theta\circ \Omega(({x}''_j)^\dagger)  = ({x}''_j)^\dagger$ and $\Theta\circ \Omega((\overline{x''_j})^\dagger)  = (\overline{x''_j})^\dagger$.

Since $\Theta\circ \Omega$ acts as the identity on the generators of $R_{\frak{F}_{Y_1}}$,
 it acts as  the identity map on the algebra $R_{\frak{F}_{Y_1}}$. 
Therefore,
 the maps $\Theta: \widehat{\widetilde{R}}_{W,S} \to R_{\frak{F}_{Y_1}}$ 
and $\Omega : R_{\frak{F}_{Y_1}} \to \widehat{\widetilde{R}}_{W,S}$ are inverse homomorphisms,
and we have unital algebra isomorphisms:
$$\widehat{\widetilde{R}}_{W,S} \cong  R_{\frak{F}_{Y_1}} \cong M_{n + m+2}({L}_a) \bigoplus M_{n + 1}(\Bbb{Q}) \bigoplus M_{m + 1}(\Bbb{Q}).$$
\end{proof}

\begin{corollary}  Let $(W, S)$ be a Coxeter system.
We have the following unital ring isomomorphisms:
$$\widehat{\widetilde{R}}_{W,S} \cong \begin{cases} 
M_2({L}_a) \oplus \Bbb{Q} \oplus \Bbb{Q} & \mbox{if $(W, S)$ is of type $I_2(a), a  > 6$, where $a$ is even,}\\
M_p(\Bbb{Q}) \oplus M_{p-1}(\Bbb{Q}) \oplus \Bbb{Q} & \mbox{if $(W, S)$ is of type $B_p$, $p \geq 2$,}\\
M_4(\Bbb{Q}) \oplus M_{2}(\Bbb{Q}) \oplus M_{2}(\Bbb{Q}) & \mbox{if $(W, S)$ is of type $F_4$,}\\
\end{cases}$$
 (where $\displaystyle{L}_n \cong \bigoplus_{N\in\mathbb{N}_{\geq 3},N\mid n}\Bbb{Q}_N, \ n \geq 3$ is a direct sum of field extensions of $\Bbb{Q}$).
\end{corollary}
\begin{proof}

\noindent{\bf Case 1:} If $(W, S)$ is a Coxeter system of type $I_2(a), a \geq 6$, $a$ even, we can apply Theorem \ref{general} with $n = 0$ and $m = 0$. This gives that 
$$ \widehat{\widetilde{R}}_{W,S} \cong M_2({L}_a) \oplus \Bbb{Q} \oplus \Bbb{Q}. $$
\noindent{\bf Case 2:} If $(W, S)$ is a Coxeter system of type $B_p, p \geq 2$,  we can apply Theorem \ref{general} with $a = 4, n = p - 2$ and $m = 0$. This gives that 
$$ \widehat{\widetilde{R}}_{W,S} \cong M_p({L}_4) \oplus M_{p - 1}(\Bbb{Q}) \oplus \Bbb{Q}  \cong M_p(\Bbb{Q}) \oplus M_{p - 1}(\Bbb{Q}) \oplus \Bbb{Q}, $$
since $L_4 \cong \Bbb{Q}$. \\
\noindent{\bf Case 3:} If $(W, S)$ is a Coxeter system of type $F_4$ we can apply Theorem \ref{general} with $a = 4, n = 1$ and $m = 1$. This gives that 
$$ \widehat{\widetilde{R}}_{W,S} \cong M_4({L}_4) \oplus M_{2}(\Bbb{Q}) \oplus M_{2}(\Bbb{Q}) \cong M_4(\Bbb{Q}) \oplus M_{2}(\Bbb{Q}) \oplus M_{2}(\Bbb{Q}). $$
\end{proof}
 We summarize the results on the finite irreducible Coxeter  groups as follows:
 \begin{theorem}\label{complete}
 Let $(W, S)$ be a Coxeter system with rank $N$, where $W$ is a finite irreducible Coxeter group.  We have the following unital ring isomomorphisms:
 $$\widehat{R}_{W,S} \cong \begin{cases} M_N(\Bbb{Q}) & \mbox{if $(W, S)$ is of type $A_n (n \geq 1), B_n (n \geq 2), D_n (n \geq 4), E_6, E_7, E_8, F_4$ or $G_2$,}\\
M_N(\Bbb{Q}_5) & \mbox{if $(W, S)$ is of type $H_3,$ or $H_4$,}\\
M_2(\Bbb{Q}_n) & \mbox{if $(W, S)$ is of type $I_2(n), n \geq 5, n \not= 6$,}\\
\end{cases}$$
(where $\Bbb{Q}_n$ if the field extension of $\Bbb{Q}$ obtained by adjoining $4\cos^2\left(\frac{\pi}{n}\right)$) \\
and
 $$\widehat{\widetilde{R}}_{W,S} \cong \begin{cases} M_N(\Bbb{Q}) & \mbox{if $(W, S)$ is of type $A_n (n \geq 1), D_n (n \geq 4), E_6, E_7$ or $E_8$,}\\
M_N(\Bbb{Q}_5) & \mbox{if $(W, S)$ is of type $H_3,$ or $H_4$,}\\
M_2({L}_n) & \mbox{if $(W, S)$ is of type $I_2(n), n \geq 5$, where $n$ is odd,}\\
M_2({L}_n) \oplus \Bbb{Q} \oplus \Bbb{Q} & \mbox{if $(W, S)$ is of type $I_2(n), n \geq 6$, where $n$ is even,}\\
M_n(\Bbb{Q}) \oplus M_{n-1}(\Bbb{Q}) \oplus \Bbb{Q} & \mbox{if $(W, S)$ is of type $B_n$, $n \geq 2$,}\\
M_4(\Bbb{Q}) \oplus M_{2}(\Bbb{Q}) \oplus M_{2}(\Bbb{Q}) & \mbox{if $(W, S)$ is of type $F_4$,}\\
\end{cases}$$
 (where ${L}_n \cong \bigoplus_{N\in\mathbb{N}_{\geq 3},N\mid n}\Bbb{Q}_N, \ n \geq 3$ is a direct sum of field extensions of $\Bbb{Q}$).
 \end{theorem}
\vskip .1in
\noindent
\centerline{DECLARATIONS}
\vskip .1in
\noindent
{\bf Ethical Approval} Not Applicable.
\vskip .1in
\noindent
{\bf Funding} Not Applicable.
\vskip .1in
\noindent
{\bf Availability of Data and Materials} Not Applicable.

\end{document}